\documentclass[12pt]{amsart}

\usepackage{amsmath,amsthm,amscd,amsfonts,amssymb,epic,eepic,bbm,tikz-cd,mathtools}
\usepackage[pagebackref,colorlinks=true,linkcolor=blue,citecolor=blue]{hyperref}
\usepackage[noabbrev]{cleveref}

\allowdisplaybreaks
\newtheorem{thm}{Theorem}[section]
\newtheorem{cor}[thm]{Corollary}
\newtheorem{conj}[thm]{Conjecture}
\newtheorem{lem}[thm]{Lemma}

\newtheorem{prop}[thm]{Proposition}
\theoremstyle{remark}

\newcounter{remarkscounter}

\numberwithin{equation}{section}

\newcommand{\Sp}{\mathrm{Sp}}

\newcommand{\N}{\mathbb{N}}

\newcommand{\Z}{\mathbb{Z}}
\newcommand{\Q}{\mathbb{Q}}
\newcommand{\R}{\mathbb{R}}

\newcommand{\C}{\mathbb{C}}

\newcommand{\Sh}{\operatorname{Sh}}

\let\clos\overline

\let\emptyset\varnothing

\let\epsilon\varepsilon

\let\phi\varphi

\newcommand*{\mb}[1]{\ensuremath{\mathbb{#1}}}
\newcommand*{\bm}[1]{\ensuremath{\boldsymbol{#1}}}
\newcommand*{\mc}[1]{\ensuremath{\mathcal{#1}}}

\newcommand*{\bs}{\backslash}

\DeclareMathOperator{\ad}{ad}

\DeclareMathOperator{\End}{End}

\DeclareMathOperator{\Gal}{Gal}

\DeclareMathOperator{\Lie}{Lie}

\DeclareMathOperator{\trdeg}{tr.deg.}

\DeclarePairedDelimiter{\norm}{\lVert}{\rVert}
\DeclarePairedDelimiter{\abs}{\lvert}{\rvert}

\DeclarePairedDelimiter{\set}{\{}{\}}

\newcommand{\quash}[1]{}

\theoremstyle{definition}
\newtheorem{defn}[thm]{Definition}

\renewcommand{\bar}{\overline}

\numberwithin{equation}{subsection}

\renewcommand{\hat}{\widehat}

\allowdisplaybreaks

\begin{document}

\title{Algebraic Independence of Special Points on Shimura Varieties}
\author{Yu Fu and Roy Zhao}
\subjclass[2020]{Primary: 11G18; Secondary: 03C64, 11G10, 14G35}

\begin{abstract}
Given a correspondence $V$ between a connected Shimura variety $S$ and a commutative connected algebraic group $G$, and $n \in \N$, we prove that the $V$-images of any $n$ special points on $S$ outside a proper Zariski closed subset are algebraically independent. Our result unifies previous unlikely intersection results on multiplicative independence and linear independence. We prove multiplicative independence of differences of singular moduli, generalizing previous results by Pila--Tsimerman, and Aslanyan--Eterovi\'c--Fowler. We also give an application to abelian varieties by proving that the special points of $S$ whose $V$-images lie in a finite-rank subgroup of $G$ are contained in a finite union of proper special subvarieties of $S$, only dependent on the rank of the subgroup. In this way, our result is a generalization of the works of Pila--Tsimerman and Buium--Poonen.
\end{abstract}

\maketitle

\section{Introduction}

\subsection{Statement of the results}
Let $S$ be a connected Shimura variety, $G$ be a commutative connected algebraic group. Let $V \subset S \times G$ be a proper irreducible closed subvariety such that the projection on each factor is dominant, and $V$ is finite over $S$. Given a point $(s, g) \in V \subset S \times G$, we call $g$ a \emph{$V$-image of $s$}. Moreover, we say that a series of points $g_1, \dots, g_n \in G$ are algebraically independent if they do not all lie in a proper algebraic subgroup.

The special subvarieties in $G$ are torsion cosets of connected algebraic subgroups and the special subvarieties in a Shimura variety are described in Section 2.A. We prove the independence of the $V$-images of special points of $S$. Previous work in the literature proved the linear independence of special points when $V$ is a correspondence from a modular curve or a Shimura curve to an elliptic curve, i.e. the linear dependences among Heegner points, which play an important role in the study of the Birch and Swinnerton-Dyer Conjecture. See \cite{BP}, \cite{Bal}, \cite{PT22} for results in this direction. There have also been results proving multiplicative independence of special points when $V$ is the graph of the $j$-function from $Y(1) \to \mb{G}_m$ (see \cite{PT17}, \cite{AEF23}). The aim of this article is to provide a unified framework to generalize the previous works to the study of algebraic independence of special points in general commutative groups. 

To state our main theorem, we give the definition of exemplary components following \cite{PT22}. 

\begin{defn}
	 With the notation as above, let $\pi_{S}$ and $\pi_{G}$ be the projections of $S \times G$ onto the first and second factors respectively. 
  \begin{itemize}
      \item A \emph{distinguished component} is an irreducible component of $W \subset V \cap (S' \times G')$, where $S' \times G'$ is a special subvariety of $S \times G$, such that $\pi_{S}(W) = S'$.
	
	\item Let $W$ be a distinguished component and let $B \supset \pi_{G}(W)$ be the smallest special subvariety of $G$ containing $\pi_{G}(W)$. We say that $W$ is \emph{exemplary} if there are no larger distinguished components $W' \supset W$ such that $\pi_{G}(W') \subset B$.

 \item We say a distinguished component $W$ in $G$ is \emph{dependent} if $G' \subset G$ is strict. We note that the unique non-dependent exemplary component is $V$ itself.
  \end{itemize}
  
\end{defn}

Pila and Tsimerman \cite[Theorem 1.1]{PT22} were able to deal with the case where $S=Y(1)^n$ is the $n$-copy of a modular (or Shimura) curve and $G=E^n$ is the self-$n$-product of an elliptic curve for $n\ge 1$. We prove the following theorem, which is a generalization of their main theorem to arbitrary (connected) Shimura varieties and commutative groups. The theorem describes all algebraic dependence among $V$-images of $n$ special points for any $n \ge 1$.

\begin{thm}\label{thm:main}
	Let $S$ be a connected Shimura variety, $G$ be a commutative connected algebraic group. Suppose $V \subset S\times G$ is an irreducible subvariety that maps finitely to $S$. Then, there are only finitely many exemplary components in $V$.
\end{thm}
 
Our proof of Theorem \ref{thm:main} follows a similar path to \cite{PT22}, via point-counting on definable sets in an o-minimal structures and using an Ax--Schanuel theorem in a suitable form. The point-counting result follows from the work of Habegger and Pila in \cite{HP16} and the Ax--Schanuel theorem follows from Bl\'azquez-Sanz, Casale, Freitag, and Nagloo in \cite{Bla23}.

Results in this direction have various kinds of applications. As a first consequence of Theorem \ref{thm:main}, we prove the following corollary and a weaker version, Corollary \ref{cor:LookAt2}, which might be easier to apply in practice.

\begin{cor}\label{indpt1}
    Let $S, G, V$ be as above. Suppose moreover that the projections of $V$ to $S$ and $G$ are both dominant. Fix $n \in \N$. There exists a proper Zariski closed subset $S' \subset S^n$ such that for any $n$ points $s_1, \dots, s_n \in S$ satisfying $(s_1, \dots, s_n) \not \in S'$, any $V$-images $g_1, \dots, g_n$ of $s_1, \dots, s_n$ are algebraically independent.
\end{cor}

\begin{proof}
    For $n \ge 1$, let $\bm{s}=(s_1, \cdots, s_n)$ be a special point in $S^n$ and let $\bm{g}=(g_1, \cdots, g_n)$ be a point in $G^n$ such that $(\bm{s}, \bm{g}) \in V^n$. It follows from the definition that if $(\bm{s}, \bm{g}) \in W$ for some dependent exemplary component then $g_1, \cdots, g_n$ are algebraically dependent in $G$. Conversely, if $g_1, \cdots, g_n$ are algebraically dependent in $G$, then $(\bm{s}, \bm{g})$ is a distinguished point and hence contained in some dependent exemplary component. 
\end{proof}

Let $D$ be a positive integer. For a set $s_1, \cdots, s_n$ of special points in $S$, we introduce the notion of $D$-independent in Definition \ref{D independent}. The discriminant $\Delta(s)$ of a special point is defined in Definition \ref{discriminant}. We prove several results towards linear independence in abelian varieties of dimension $g \le 3$, utilizing the classification of special subvarieties in the moduli space of abelian varieties $\mathcal{A}_g$.

\begin{cor}\label{linearly independence of CM points}
Suppose that $S \subset \mc{A}_g$ for $g \le 3$. For $n \ge 1$ there exists a positive integer $D = D(n,S,G,V)$ such that if $s_1, \dots, s_n$ are $D$-independent, then any $V$-image $x_1, \dots, x_n$ of $s_1, \dots, s_n$ are algebraically independent in $G$.
\end{cor}

Taking $V$ to be the graph of the function $f \colon Y(1)^n \times Y(1)^{n} \to \C^n$ given by
\[f(z_1, \dots, z_n, w_1, \dots, w_n) = (z_1 - w_1, \dots, z_n - w_n),\]
we get the following generalization of \cite[Thm. 1.2]{PT17} and \cite[Thm. 1.1]{AEF23}.

\begin{cor}
    For $n \ge 1$, there exists a positive integer $D = D(n)$ such that if $x_1, \dots, x_n,$ $y_1, \dots, y_n \in Y(1)$ are singular moduli that are $D$-independent, then there does not exist $(a_1, \dots, a_n) \in \Z^n \bs \set{\bm{0}}$ such that
    \[\prod_{i = 1}^n (x_i - y_i)^{a_i} \in \mu_\infty,\]
    where $\mu_\infty$ is the set of roots of unity.
\end{cor}

 We prove also that on low dimensional abelian varieties, the intersection of any finite-rank subgroup of $G$ with the set of $V$-images of special points on $S$ is contained in a finite union of proper special subvarieties. This generalizes Pila and Tsimerman's result \cite[Corollary 1.4]{PT22} to higher dimensional abelian varieties, and Buium and Poonen's \cite[Theorem 1.1]{BP} in a uniform way, depending only on the rank of the subgroup.  

\begin{thm}\label{Andre-Oort-Mordell-Lang}
    Suppose $G$ is defined over a field $K$ of characteristic $0$ and suppose the Shimura datum for $S$ has $G_S$ simple. Let $\Gamma \subset G(K)$ be a subgroup of rank $r$ and let $\Gamma'$ be its division group. Then, there exist $N(r, S, G, V)$ proper special subvarieties of $S$ such that the set of special points of $S$ with a $V$-image in $\Gamma'$ is contained in a union of these proper special subvarieties.
\end{thm}

A similar argument to \cite[Proposition 3.2]{PT22} proves that Theorem \ref{thm:main} is a consequence of the Zilber--Pink Conjecture. Also, Theorem \ref{Andre-Oort-Mordell-Lang} may be thought as a weak version of the Andr\'e--Oort--Mordell--Lang conjecture in the case of subgroups of finite rank formulated by Baldi in \cite[Conjecture 5.2]{Bal}. The original formulation of the Andr\'e--Oort--Mordell--Lang conjecture is due to Pink, which can be found in \cite{Pink05}.

\subsection{Organization of the paper}
The notions of Shimura varieties, special and weakly special subvarieties are explained in Section 2. We also prove a mixed version of the Ax--Schanuel theorem for the uniformization of the product $S\times G$ in the form needed for our main theorem in this section. In Section 3, we collect and prove some arithmetic estimates. The main theorem is proved in Section 4 and its applications and corollaries are proved in Section 5. 

\section{Shimura varieties and the mixed Ax--Schanuel}
Let $S$ be an algebraic variety and $q \colon X \to S$ be a universal covering of $S$. The Ax--Schanuel theorem gives information about the bi-algebraic varieties, algebraic varieties $V \subset X$ such that $q(V)$ is also algebraic, for the transcendental covering map $q$. The rough statement of the Ax--Schanuel theorem is that varieties of $X \times S$ should have a proper intersection with the graph of $q$ unless the projection of $V$ to $S$ is contained within a bi-algebraic variety. These bi-algebraic varieties are precisely the weakly special subvarieties. We first give precise definitions of bi-algebraic varieties and results in the case when $S$ is a commutative algebraic group or a Shimura variety, and then prove the result when $S$ is a product of the two.

\subsection{Shimura Varieties}\label{Section 2.1}
Let $S:=\Sh_K(G_S, X)$ be the (connected) Shimura variety associated with a connected Shimura datum $(G_S, X)$, where $G$ is a semisimple group of adjoint type and $K$ is a compact open subgroup of $G(\mathbb{A}_f).$ More information about Shimura varieties can be found in Milne (\cite{Mil05}).

For any Shimura subvariety $Z$ of $S$ and any $a \in G(\mathbb{A}_f)$, we refer to any irreducible component of the Hecke correspondence $T_{K,a}(Z)$ as a \textit{special subvariety} of $S$. A \textit{special point} is a special subvariety of dimension zero. 

 For any $x \in X$, let $\mathbf{M}:= \operatorname{MT}(x)$ be the Mumford-Tate group of $x$, which is defined as the smallest $\Q$-subgroup $H$ of $G$ such that $x$ factors through $H_{\mathbb{R}}$. Let $X_{\mathbf{M}}$ be the $\mathbf{M}(\mathbb{R})$-conjugacy class of $x$. Then the image of $X_\mathbf{M}$ in $S = \Gamma_g \backslash X $ is a special subvariety of $S$. It is not hard to see that every special subvariety of $S$ arises this way. The action of $\mathbf{M}(\mathbb{R})$ on $X_\mathbf{M}$ factors through the adjoint group $\mathbf{M}^{\operatorname{ad}}(\mathbb{R})$, which is a direct product of its $\Q$-simple factors. Therefore, we can write $X_\mathbf{M}$ as a product
$$X_{\mathbf{M}}=X_1 \times X_2,$$
corresponding to the action of $M^{\operatorname{ad}}(\mathbb{R})$. A \textit{weakly special subvariety} is the image of the fiber $\{x_1\} \times X_2$ or $X_1 \times \{x_2\}$ for any $x_1 \in X_1$ or $x_2 \in X_2$. Therefore, a weakly special subvariety of $S$ is a special subvariety if and only if it contains a special point. The weakly special subvarieties of $S$ are precisely those subvarieties that are totally geodesic in $S$ (\cite[Section 4]{Mo98}), and they are also precisely the bi-algebraic subvarieties under the uniformization map $q \colon X \to S$. More detailed information on weakly special subvarieties of Shimura varieties can be found in \cite{Ull11}.

The Ax--Schanuel theorem for Shimura varieties is a statement about the functional transcendence of the uniformization map $q$.

\begin{thm}[{\cite[Theorem 1.1]{Mok19}}]
	Let $V \subset X \times S$ be an algebraic subvariety and let $D \subset X \times S$ be the graph of $q \colon X \to S$. Let $U \subset V \cap D$ be an irreducible component such that
	\[\dim U > \dim V - \dim S.\]
	Then, the projection of $U$ to $S$ is contained in a proper weakly special subvariety of $S$.
\end{thm}

\subsection{Commutative Algebraic Groups}\label{Section 2.2}
Let $G$ be a connected commutative algebraic group of dimension $g$ defined over $\C$ and let $\exp \colon \C^g \to G$ be the exponential map from $\Lie G \cong \C^g$ to $G$.

\begin{defn}\label{special subvarieties abelian}
	For each connected subgroup $B \subset G$ and closed point $p \in G(\C)$, we say that the algebraic subvariety $B + p \subset G$ is a \emph{weakly special subvariety of $G$}. We say that it is a \emph{special subvariety of $G$} if $p \in G_{\mathrm{tors}}$ is a torsion point.
\end{defn}

The special points of $G$ are special subvarieties of dimension $0$, and hence they are precisely the torsion points of $G$. As before, a weakly special subvariety of $G$ is special precisely when it contains a special point.

The Ax--Schanuel theorem for commutative groups is a statement about covering map given by the $\Lie$-exponential $\exp$.

\begin{thm}[{\cite[Theorem 3]{Ax72}}]
	Let $V \subset \C^g\times G$ be an algebraic subvariety and let $D \subset \C^g\times G$ be the graph of $\exp \colon \C^g \to G$. Let $U \subset V \cap D$ be an irreducible component such that
	\[\dim U > \dim V - \dim G.\]
	Then the projection of $U$ to $G$ is contained in a proper weakly special subvariety of $G$.
\end{thm}

\subsection{Mixed Ax--Schanuel}
We will prove a version of the Ax--Schanuel theorem for the product uniformization map $q \times \exp \colon X \times \C^g \to S\times G$. This was proven by Pila and Tsimerman in \cite{PT22} for the case when $S$ is a product of modular curves and $G$ is a product of elliptic curves. In that case, it was shown that the Ax--Schanuel theorem for the product follows from an Ax--Schanuel theorem for mixed Shimura varieties by Gao in \cite{Gao20}. Using the same method, one can prove that the Ax--Schanuel theorem holds for the product $S\times G$ whenever $S$ is a Shimura variety of abelian type and $G$ is an abelian variety. However, to prove it for general Shimura varieties and general commutative groups, we will need stronger machinery from Bl\'azquez-Sanz, Casale, Freitag, and Nagloo (\cite{Bla23}). The translation into the language of (weakly) special subvarieties is given by Chiu (\cite{Chi22}).

\begin{defn}
	Let $S$ be a Shimura variety and $G$ a connected commutative algebraic group. We say that an algebraic subvariety $V \subset S\times G$ is a \emph{weakly special subvariety of $S\times G$} if there exist weakly special subvarieties $S_H \subset S$ and $B \subset G$ such that $V = S_H \times B$. We say that $V$ is a \emph{special subvariety} if $S_H, B$ are both special.
\end{defn}

\begin{thm}\label{thm:ax-schanuel}
	Let $V \subset (X \times \C^g) \times (S\times G)$ be an algebraic subvariety and let $D \subset (X \times \C^g) \times (S\times G)$ be the graph of $q \times \exp \colon (X \times \C^g) \to (S\times G)$. Let $U \subset V \cap D$ be an irreducible component such that
	\[\dim U > \dim V - \dim (S\times G).\]
	Then the projection of $U$ to $S\times G$ is contained in a proper weakly special subvariety of $S\times G$.
\end{thm}
\begin{proof}
    Let the Shimura datum of $S$ be $(G_S, X)$. We may take $S$ to be a connected Shimura variety and $G$ to be a derived group. Applying \cite[Thm. 3.6]{Bla23} to $D \subset (X \times \C^g) \times (S\times G)$ gives that the projection of $U$ in $S\times G$ is contained in a proper subvariety whose Galois group is a strict algebraic subgroup $H \subset G_S(\C) \times \mb{G}_a(\C)^g$. Since $G$ is a derived group, it has no abelian quotients and hence Goursat's lemma says that the projection of $H$ to $G(\C)$ or $\mb{G}_a(\C)^g$ is not onto. By \cite[Thm. 3.2]{Chi22} for the Shimura variety side and \cite[Thm. 3]{Ax72} for the commutative group side, that means the projection of $U$ to $S$ or $G$ is contained in a proper weakly special subvariety.
\end{proof}


\section{Bounds on the Galois orbits}

In this section we give some arithmetic estimates which will be used later in the proof of the theorems.

Fix a choice of fundamental domain $\mathcal{F}_S$ for the uniformization $q : X \to S$ of the Shimura variety $S$. Let $F$ be the number field over which $S$ admits a canonical model. The degree of $F$ is bounded in terms of the datum $(G_S,X)$ and $K$. All the special points of $S$ are algebraic points defined over abelian extensions of $F$.

 For a special point $s \in S$, let $x$ be a preimage of $s$, i.e. $x$ is a pre-special point. By definition, the Mumford--Tate group of a special point is an algebraic torus. Let $K^{m}_{\mathbf{M}}$ be the maximal compact open subgroup of $\mathbf{M}(\mathbb{A}_f)$ and $K_\mathbf{M}$ the compact open subgroup $K \cap \mathbf{M}(\mathbb{A}_f)$ of $\mathbf{M}(\mathbb{A}_f)$. Let $E$ be the splitting field of $\mathbf{M}$. Since $G_S$ is of adjoint type, $E$ is a (Galois) CM field. Let $D_E$ be the absolute value of the discriminant of $E$.
\begin{defn}\label{discriminant}

   The discriminant of $s$ is $$\Delta(s) \coloneqq [K^{m}_{\mathbf{M}}: K_\mathbf{M}]D_{E}.$$
   \end{defn}
   
We want to estimate the heights and degrees of special and pre-special points in terms of $\Delta(s).$ Fortunately, recent progress allows us to have the following proposition.
 
   \begin{prop}\label{arithmetic estimates}
 	Let $s \in S$ be a special point with discriminant $\Delta(s)$ and  let $x$ be a preimage of $s$ under the uniformization map $q$. Let $h$ be the canonical height on $S$, which is a logarithmic Weil height, and let $H$ be the multiplicative Weil height for a fixed realization of $X$. 
 	Then we have 
 	\begin{itemize}
 		\item[(a)] $h(s) \le |\Delta(s)|^{o(1)}$;
 		\item[(b)] $H(x) \le C_1|\Delta(s)|^{C_2}$;
 		\item[(c)] $[F(s): F] \ll |\Delta(s)|^{1/2+ \epsilon}$ for any $\epsilon > 0$;
 		\item[(d)] $[F(s): F] \gg C_3 |\Delta(s)|^c$ for some fixed $c>0$.
 	\end{itemize}
 	Where the constants $C_1, C_2, \dots$ are depend on $G_S, X, F, \mathcal{F}_S$ and the realization. 
 \end{prop}

\begin{proof}
\begin{itemize}
	\item[(a)] Follows from \cite[Theorem 9.11]{PST};
 		\item[(b)] See \cite[Theorem 1.1, 4.1]{DO};
 		\item[(c)] Follows from the Brauer-Siegel theorem for arithmetic tori in \cite[Theorem 1.3]{Tsi12};
 		\item[(d)] See \cite[Theorem 1]{BSY}.
\end{itemize}
	
\end{proof}

Over $\clos{\Q}$, any commutative connected group $G$ can be written as a product of a semi-abelian variety with affine space $G' \times \mb{G}_a^m$ (see \cite[Prop. 5.1.12]{JJ14}). We will need height bounds on the semi-abelian factor, so suppose $G$ is a semi-abelian variety over a number field $L$ with toric part $T$ and abelian quotient $\pi : G \to A$. Let $g$ be the dimension of $A$. Let $K$ be any number field. Note that the Weil height on a semi-abelian variety might be negative. However, the canonical height $\hat{h}_{L}$ on $G$ defined by K\"uhne \cite[Sec. 3]{Kuh}, with respect to a $T$-effective line bundle $(M,\varrho)$ on a $T$-equivariant compactification $\bar{T}$ and an ample line bundle $N$ on $A$, remedies this issue. As for abelian varieties, the zero set of $\hat{h}_{L}$ coincides with the torsion points of $G$. Let
$$\eta:=\eta(G,L)=\operatorname{inf} \hat{h}_{L}(P),$$
where the infimum is taken over all non-torsion $P$ in $G(K)$, and we write
$$\omega := \omega(G, K)$$
for the cardinality of the torsion group of $G(K)$. Note that $\eta$ and $\omega$ depend on the embedding of $G$ in a projective space. Suppressing this dependence in our notation, we can prove the following theorem using an argument similar to the proof of \cite[Theorem $\mb{G}_m, A$]{Mas88}.
 \begin{thm}
 Suppose $P_1, \ldots, P_n$ on $G(K)$ have canonical heights at most $q \geq \eta$. Then the relation group of $P_1, \ldots, P_n$ generated by $$\bm{m}: m_1P_1 + \cdots m_nP_n= 0_G$$ satisfy
$$
|m_i| \leq n^{n-1} \omega(q / \eta)^{n-1} .
$$
\end{thm}

\begin{proof}
    We refer to \cite[Lem. 8]{Kuh} for the definition and properties of the canonical height on a semi-abelian variety. By definition, $$\hat{h}_L=\hat{h}_{G(M, \varrho)}+\hat{h}_{\pi^* N},$$ where $\hat{h}_{G(M, \varrho)}$ is linear and $\hat{h}_{\pi^* N}$ is quadratic. We define a convex distance function $f$ on $\Z^{n}$ by
    $$f(\bm{m})= \hat{h}_{G(M, \varrho)}(\bm{m})+ \hat{h}_{\pi^* N} (\bm{m})^{1 / 2}.$$
    Let $\Gamma=\Gamma(f)$ be the set consisting of all $\bm{m}$ such that $m_1 P_1+\ldots+m_n P_n$ is a torsion point of $G(K)$, and take $\Gamma_0$ as the relation group of $P_1, \ldots, P_n$. We can choose $E=q, \epsilon=\eta^{1 / 2}$, and the theorem follows immediately from \cite[Prop.]{Mas88}.
    
\end{proof}

General semi-abelian subvarieties of $G$ are generated by relation groups where the $m_i \in \End(G)$. We accommodate for $\End(G) \ne \Z$ in the following way. Let $\set{1, \tau_1, \dots, \tau_k} \in \End(G)$ be a set of generators. Then, we can work in $G^{(k + 1)n}$ with $P_i, \tau_{1}P_i, \cdots \tau_{k}P_i$ for $1 \le i \le n$. Due to the dominance of a height associated to an ample divisor (see \cite[Chap. 4, Prop. 5.4]{Lan83} or \cite[Lemma 6]{Poo00}), 
\begin{align}
	\hat{h}_L(\tau_i(P)) \le C^{i}\hat{h}_L(P)
\end{align}
where each $C^{i}$ is an absolute constant depending on $\tau_i$ and the embedding of $G$.  Let $$\norm*{a_1 +\sum_{j=1}^{k}a_{j + 1}\tau_j}= \textup{max}_{1 \le j \le k + 1}\{|a_j|\}$$ denote the norm of $\bm{a}=a_1 +\sum_{j=1}^{k}a_{j + 1}\tau_j$.
\begin{cor}\label{bound linear relations}
\begin{align}\label{upperbound m_i}
	\norm*{m_i} \le \bar{C} (2ng)^{2ng-1}\omega(q/\eta)^{(2ng-1)/2}
\end{align}
where $\bar{C}$ is taken as the maximal of $(C^{i})^{(2ng-1)/2}$ and $q:=\textup{max}_{1 \le i \le n}\{\hat{h}_L(P_i)\}$.
\end{cor}

To estimate \ref{upperbound m_i} in terms of degree of the algebraic points, we need some estimates for $\eta, \omega$, which is given for tori and abelian varieties in \cite{Mas88}. Let $D=[K:\Q]$ and $\mc{L}=\log (D+2)$ we have:
$$\eta \geq C^{-1} D^{-(2 g+1)} \mathcal{L}^{-2 g}$$
and
$$\omega \leq C D^g \mathcal{L}^g.$$

To bound the height of rational points on a suitable definable set in Section 4, we need the following upper bound of the norm of the generating set for the linear relations satisfied by $n$ fixed $V$-images of special points in terms of their discriminants.     
\begin{prop}\label{prop:bound endomorphism size}
	Let $(s_1,x_1), \cdots, (s_n, x_n) \in S\times G \times \mb{G}_a^m$ be points of $V$ with $s_i$ special and with discriminants $\Delta(s_i)$. Define the complexity of $s=(s_1, \dots, s_n)$ by $$\bm{\Delta}(s)=\bm{\Delta}(s_1, \cdots, s_n)= \operatorname{max}|\Delta(s_i)|.$$ Then there are constants $C, C^{\prime}, c^{\prime}$ depending on $S,G,V,n$ such that for $\bm{\Delta} \ge C$, there is a generating set for the linear relations satisfied by the $x_i \in G$ such that $$||m_i|| \le C^{\prime} \bm{\Delta}(s)^{c^{\prime}}.$$
\end{prop}
\begin{proof}
Since $V$ is a correspondence in $S\times G$, $(s,x) \in V$ is an algebraic point, using standard properties of heights we have that	
$$H(x) \le C_{1}^{\prime}H(s)^{c}$$ and since the projection of $V$ to $S$ is finite and both projections are dominant on each factor, we have
$$ [L(x):L] \le C_{2}^{\prime}[F(s):F].$$
Take $C^{\prime}=\operatorname{max}\{C_{1}^{\prime}, C_{2}^{\prime}\}.$ Therefore the degree $D=[L^{\prime}: \Q]$ of the field of definition of $x_1, \cdots , x_i$ is bounded in terms of $\bm{\Delta}$ by Proposition \ref{arithmetic estimates}(c). 

By \cite[Lem. 8(a)]{Kuh}, after fixing the embedding corresponding to the $T$-effective and ample line bundles on $G$, the differences of $\hat{h}_{L}(P)$ and $h(P)$ are bounded globally on $G(\overline{\Q})$ by an absolute constant $\delta$. Therefore for $h:=\textup{max}\{h(x_1), \cdots, h(x_n)\}$ sufficiently large, we have $$\eta \le h - \delta$$ and $$h \ge \delta$$ which implies $$q \le 2h.$$ By Proposition \ref{arithmetic estimates}(a) $\eta$ and $q $ are bounded in terms of $\bm{\Delta}$, and the proposition follows from Corollary \ref{bound linear relations}.

\end{proof}

\section{Exemplary Components}

\subsection{Proof of the main theorem over \texorpdfstring{$\bar{\Q}$}{Qbar}}
Now we can prove Theorem \ref{thm:main}. We first prove the theorem for varieties $V$ defined over $\clos{\Q}$, and then show the result for all $V$. 

First, we prove a criterion for inclusion of complex algebraic varieties defined over $\overline{\Q}$ and a finiteness result for pre-special subvarieties.

\begin{lem}\label{lem:trascendentalPoints}
    Let $V \subset \C^n \times \C^m$ and $W \subset \C^n$ be irreducible algebraic varieties defined over $\clos{\Q}$. There exists a finite set of points $P_1, \dots, P_{m + 1} \in \C^n$ such that for any $Q \in \C^m$, we have $W \times \set{Q} \subset V$ if and only if $P_i \times Q \in V$ for all $i$.
\end{lem}
\begin{proof}
    Let $d = \dim W$. We can find a generic point $P_1 \times P_2 \times \cdots \times P_{m + 1} \in W^{m + 1}$ such that $\trdeg_{\Q} \Q(P_1, \dots, P_{m + 1}) = d(m + 1)$, and we claim that these $P_i$ satisfy the above property. The forward direction of the if and only if is clear. Now suppose that $P_i \times Q \in V$ for all $i$. Viewing $P_i$ as a set, let $\Tilde{P_i} \subset P_i$ be a minimal transcendental basis for $\trdeg_\Q \Q(P_i)$. Write $Q = (q_1, \dots, q_n) \in \C^n$ and let $Q_i \subset \set{1, 2, \dots, n}$ be the set of indices $j$ so that $q_j$ is algebraic over $\Q(\Tilde{P_i})$. By construction, $\abs{\Tilde{P_i}} = d$ and the $\Tilde{P_i}$ are algebraically independent. This implies that $Q_i \cap Q_j = \emptyset$ for $i \neq j$ and so $\sum_{i = 1}^{m + 1} \abs{Q_i} \le m$. Thus, there exists an $i$ such that $\abs{Q_i} = 0$ and for this $i$ we have
    \[\trdeg_\Q \Q(P_i, Q) = \trdeg_\Q \Q(P_i) + \trdeg_\Q \Q(Q).\]
    Let $X = \clos{\set{P_i}}^{\mathrm{Zar}}$ and $Y = \clos{\set{Q}}^{\mathrm{Zar}}$ be the $\clos{\Q}$-Zariski closure of $P_i$ and $Q$. Then $\dim X = \trdeg_\Q \Q(P_i) = d$ and $\dim Y = \trdeg_\Q \Q(Q)$. Note that the Zariski closure of $P_i \times Q$ has dimension $\trdeg_\Q \Q(P_i, Q) = \dim X + \dim Y$, showing that the $\clos{\Q}$-Zariski closure of $P_i \times Q$ is $X\times Y$. But $X\times Y \subset V$ by definition and $X = W$ showing inclusion in the other direction.
\end{proof}

\begin{lem}\label{lem:finitely many families}
    Suppose $V \subset S\times G$ is an algebraic subvariety of the product of a Shimura variety with a connected commutative algebraic group. There exists a finite set $\Sigma$ of sub-Shimura datum $(H, X_H)$ of $(G_S, X)$ and splittings $(H^{\ad}, X_H^{\ad}) = (H_1, X_1) \times (H_2, X_2)$ such that if $W \subset V$ is exemplary, then there exist one such splitting in $\Sigma$ such that $\pi_S(W)$ is the image of $X_1 \times \set{x_2}$, for some $x_2 \in X_2$, under the uniformization map of $S$.

    Moreover, if $\mc{V} \to W$ is a definable family of $V \subset S\times G$, then the result still holds.
\end{lem}
\begin{proof}
    The result for a fixed $V \subset S\times G$ is proven by \cite[Prop. 3.4]{PT22} and \cite[Prop. 6.10]{DR18}. The version for family is given by the uniform Ax--Schanuel theorem, which is implied by the Ax--Schanuel theorem. This is proven in \cite[Prop. 2.20]{ES23}.
\end{proof}

We will need the following stronger form of the Pila--Wilkie point-counting theorem for families.

\begin{thm}{\cite[Cor. 7.2]{HP16}}\label{thm:point-counting}
    Let $F \subset \R^\ell \times \R^m \times \R^n$ be a definable family parametrized by the first factor $\R^\ell$. Let $\epsilon > 0$ and $k \in \N$ and let $\pi_1 \colon \R^m \times \R^n \to \R^m$ and $\pi_2 \colon \R^m \times \R^n \to \R^n$ denote the projections onto the first and second factors. There exists a constant $c = c(F, k, \epsilon) > 0$ satisfying the following.

    Let $x \in \R^\ell$ and let $F_x \subset \R^m \times \R^n$ denote the fiber of $F$ over $x$. If $T \ge 1$ and there exists a subset
    \[\Sigma \subset \set{(y, z) \in F_x : H_k(y) \le T}\]
    such that $\abs{\pi_2(\Sigma)} > cT^\epsilon$, then there exists a continuous definable function $\beta \colon [0, 1] \to F_x$ satisfying the following four properties:
    \begin{enumerate}
        \item The composition $\pi_1 \circ \beta$ is semialgebraic and its restriction to $(0, 1)$ is real analytic;
        \item The composition $\pi_2 \circ \beta$ is non-constant;
        \item $\pi_2(\beta(0)) \in \pi_2(\Sigma)$;
        \item and the restriction of $\beta$ to $(0, 1)$ is real analytic.
    \end{enumerate}
\end{thm}

\begin{thm}\label{thm:overQbar}
	Let $S$ be a connected Shimura variety and let $G$ be a connected commutative algebraic group. Suppose $V \subset S\times G$ is a proper irreducible subvariety that is finite over $S$. Moreover, suppose that $G$ and $V$ are defined over $\clos{\Q}$. Then, there are only finitely many exemplary components in $V$.
\end{thm}
\begin{proof}
	Let $\pi_S \colon V \to S$ and $\pi_G \colon V \to G$ be the projections onto the two factors. Let $W \subset V$ be an exemplary component and let $S' \coloneqq \pi_S(W)$ be the special subvariety of $S$ which $W$ maps onto, and let $G' \supset \pi_G(W)$ be the smallest special subvariety of $G$ containing $\pi_G(W)$. $S'$ is a special subvariety of $S$ and hence there is a sub-Shimura datum $(H, X_H)$ of $(G_S, X)$, a decomposition $(H^{\ad}, X_H^{\ad}) = (H_1, X_1) \times (H_2, X_2)$, and a point $y_2 \in X_2$ so that $S'$ is the image of $X_1 \times \set{y_2}$. By Lemma \ref{lem:finitely many families}, there are only finitely many choices of $X_H$ and $X_1$. Thus, to show finiteness of exemplary components, it suffices to show that for a fixed $X_1$ and fixed dimension of $G'$, there are only finitely many points $y_2 \in X_2$ arising from exemplary subvarieties.
 
    Let $\tilde{q} \colon X_1 \times X_2 \to X \to S$ denote the uniformization map of $S$ restricted to $X_1 \times X_2$ and let $q \colon F_1 \times F_2 \to S$ denote the restriction of $\tilde{q}$ to a fundamental domain $F_1 \times F_2$ of $X_1 \times X_2$. Let $\exp \colon \C^{g} \to G$ denote the uniformization map of the algebraic group and let $e \colon F_G \to G$ denote the restriction of that map to a fundamental domain. Using Lemma \ref{lem:trascendentalPoints}, let $x_1, \dots, x_r \in F_1$ denote a set of points to determine if $q(X_1 \times \set{x_2}) \times \exp(z) \subset V$. Let $\xi_1, \dots, \xi_g$ denote the coordinates in $G$ and let $\tau_1, \dots, \tau_k \in \End(G)$ be a $\Z$-basis for $\End(G)$. Take a generating set of all equations of the form
    \[\sum m_{ij} \tau_j \xi_i = 0, m_{ij} \in \Z\]
    that all the points of $G'$ satisfy, and let $G_0$ be the identity component of the algebraic subgroup of $G$ defined by these equations. We can extend each map $\tau_i \colon G \to G$ to an endomorphism of its covering space $\tilde{\tau_i} \colon \C^{g} \to \C^{g}$ satisfying $\tilde{\tau_i}(0) = 0$. Suppose that $G_0$ is cut out by $\ell$ such equations. Let
    \begin{align*}
        Y =& \left\{(y, z_1, \dots, z_g, m_{111}, \dots, m_{gk\ell}, b_1, \dots, b_\ell) \in F_2 \times F_G \times \R^{gk\ell} \times \R^\ell :\right.\\
        &\left.\forall 1 \le i \le r, (q(x_i, y), e(\bm{z})) \in V\text{ and }\forall 1 \le l \le \ell, \sum_{i, j} m_{ijl}\tilde{\tau_j}(z_i) = b_l\right\}
    \end{align*}
    and set $Z$ to be the projection of $Y$ to $F_2 \times \R^{gk\ell} \times \R^\ell$. Both $Y$ and $Z$ are definable sets. The set $Y$ parametrizes points $y \in F_2$ so that a $V$-image of $X_1 \times \set{y}$ lies within a proper special subvariety cut out by the $\bm{m}$ of $G$ (but not all choices of $\bm{m}$ and $\bm{b}$ correspond to algebraic subvarieties).

    Suppose for the sake of contradiction that there were infinitely many $\clos{\Q}$-algebraic exemplary subvarieties $W'$ with fibers over $X_1 \times X_2$, and with dimension of the smallest special subvariety of $G$ containing $\pi_G(W')$ equal to a fixed dimension $d = \dim G'$. Each one gives a $\clos{\Q}$-point $(y', \bm{m}, \bm{b}) \in Z$, and the $\Gal(\clos{\Q}/K)$-orbits of $W'$ also lie in $Z$ for $K$ the defining field of $S, G, \End(G)$. Over $\clos{\Q}$, any commutative connected group $G$ can be written as a product of a semi-abelian variety with affine space $G' \times \mb{G}_a^n$ (see \cite[Prop. 5.1.12]{JJ14}). The tuple $\bm{m}$ consists integers and for semi-abelian varieties, we may take the fundamental domain so that the real part of each $z_i$ is in the interval $[0, 1]$. Thus, each $b_l$ corresponding to an equation on the semi-abelian variety is an integer bounded by $\sum \abs{m_{ijl}}$, which by Proposition \ref{prop:bound endomorphism size}, is bounded by $C\abs{\Delta(y')}^c$. Special subvarieties of $\mb{G}_a^n$ are given by the equations $\sum m_{ij}\tilde{\tau_j}(z_i) = 0$, so $b_l = 0$ for those equations. Thus, taking the Galois orbit of $W'$ gives a point with different $y'$ as well as different $b_l$, but the same $\bm{m}$. View $Z$ is fibered over the $\bm{m}$ variable and let $\Sigma \subset Z_{\bm{m}}$ be the set of points arising from exemplary subvarieties $W'$ as well as their Galois conjugates. By Proposition \ref{arithmetic estimates}, there are at least $C' \abs{\Delta(y')}^{c'}$ points of height less than $C\abs{\Delta(y')}^{c}$. Then Theorem \ref{thm:point-counting} gives the existence of a set $R \subset Z_{\bm{m}}$ of positive dimension whose projection to $F_2$ is connected semialgebraic and whose projection to $\R^\ell$ is non-constant.

    Let $\clos{\exp} \colon \C^{g}/\tilde{G_0} \to G/G_0$ be the exponentiation map of $G/G_0$. Let $F_{G_0}$ denote the image of the fundamental domain $F_G$ under the quotient map, which will serve as a fundamental domain for $\clos{\exp}$. First, take the preimage of $R \subset Z_{\bm{m}}$ under the projection map $Y \to Z$ and then project the preimage under the map
    \[Y_{\bm{m}} \to F_2 \times F_{G_0} \times \R^\ell \to F_2 \times F_{G_0}.\]
    Let the image of $R$ under these transformations be $R' \subset F_2 \times F_{G_0}$. By applying the point counting theorem if necessary, we may take $R'$ to be connected and semialgebraic. Since the projection of $R$ to $\R^\ell$ was non-constant, the projection of $R'$ to $F_{G_0}$ is also non-constant. Moreover, we may assume that $G/G_0$ is a semi-abelian variety by expanding $G_0$ if necessary because the projection to the affine space $\mb{G}_a$ is constant ($b_l = 0$). By applying the Ax--Schanuel theorem (Theorem \ref{thm:ax-schanuel}) to the Zariski closure of $(F_1 \times R') \times (q \times \exp)(F_1 \times R')$, we get that $(q \times \exp)(F_1 \times R') \subset S\times G/G_0$ is contained in a proper weakly special subvariety. Since $R'$ contains preimages of special points, its image must lie in a proper special subvariety $S'\times G' \subset S\times G/G_0$.
    
    Let $V' \subset S\times G/G_0$ be the image of $V$ under the quotient map $G \to G/G_0$. By construction, we have that $(q \times \exp)(F_1 \times R') \subset V'$, and so $V'$ contains a proper special subvariety $S'\times G'$. However, the map $V \to S$ is finite, and hence $V' \to S$ is also finite. Hence, $G'$ is a single point, meaning that $R$ was originally of the form $R'' \times \set{g}$, with $g$ corresponding to a torsion translation of $G_0$, which by abuse of notation we will also denote $G_0$. Let $\tilde{G_0} \subset F_G$ be the preimage of $G_0$. Applying Ax--Schanuel to the Zariski closure of
    \[U \coloneqq (F_1 \times R'' \times \tilde{G_0}) \times (q(F_1 \times R'')\times G_0) \cap V,\]
    we get that $q(F_1 \times R'')\times G_0$ is contained in a proper special subvariety of $V$, which must be of the form $S'\times G_0 \subset V$. This shows that our original exemplary component $W$ was not exemplary because there is a larger Shimura variety $S'$ properly containing $\pi_S(W)$ whose $V$ image lies within $G_0 = G'$.
\end{proof}

\subsection{From \texorpdfstring{$\bar{\Q}$}{Qbar} to \texorpdfstring{$\C$}{C}}

Armed with the result over number fields, through a continuity argument, we are able to prove Theorem \ref{thm:main} over any characteristic $0$ field.

\begin{proof}[Proof of Theorem \ref{thm:main}]
    Let $F$ be the finitely generated subfield over $\Q$ which $S, G, \End(G)$, and $V$ are all defined. Then, there exists a quasi-projective geometrically irreducible variety $Z$ over a number field $K$ such that $F$ is the function field of $Z$, and $S$ is defined over $K$. By taking a Zariski open subset of $Z$, we may assume that $G$ extends to an abelian scheme $\mc{G}$ over $Z$ and $V$ extends to a variety $\mc{V}$ that is flat over $Z$. Choose a generic closed point $z_0 \in Z(\C)$ such that the field gotten by localizing at $z_0$ is isomorphic to $F$. Choose a simply connected Euclidean open neighborhood $U \ni z_0$ of $z_0$ so that the homology of $\mc{G}$ can be trivialized over $U$, and so we have $\mc{G}_U \cong G \times U$ as analytic varieties.

    Suppose for the sake of contradiction that there are infinitely many exemplary subvarieties of $V \cong \mc{V}_{z_0}$. By Lemma \ref{lem:finitely many families}, special subvarieties that arise as projections of the exemplary subvarieties of $V$ to $S$ come from finitely many splittings of special subvarieties of $S$. So, there is a splitting of sub-Shimura datum $(H_1, X_1) \times (H_2, X_2)$ with image $S_1 \times S_2 \subset S$ and infinitely many special points $p_i \in S_2$ such that there are exemplary subvarieties $W_i$ of $V$ mapping surjectively to $S_1 \times \set{p_i}$. By the Andr\'e--Oort conjecture, proven by \cite{PST}, the Zariski closure of these $S_1 \times \set{p_i}$ is a finite union $R_1, \dots, R_n \subset S$ of special subvarieties of $S$. For each $i$, let $G_i \subset G$ be the smallest special subvariety containing $\pi_G(W_i) \subset \mc{G}_{z_0}$, and we extend $G_i$ to a family $\mc{G}_i$ over $U$. For every $z \in U$, we have that the $\mc{V}_z$ image of $S_1 \times \set{p_i}$ lies within $G_{i, z}$. However, for each $z \in U(\clos{\Q})$, Theorem \ref{thm:overQbar} says that there are only finitely many exemplary components and hence the image of each $R_i$ must also lie within a proper special subvariety of $G$.

    The $\mc{V}$-images of each $R_i$ gives a family of subvarieties $\mc{R}$ of $\mc{G}$ such that for each $z \in U(\clos{\Q})$, the fiber $\mc{R}_{z}$ lies within a proper special subvariety of $\mc{G}_z$. We claim that this holds at $z_0$ as well. By replacing $\mc{R}$ with a $g$ self-sum, where $g = \dim G$, we may assume the $\mc{R}$ is a coset of an abelian subscheme of $\mc{G}$ and we still have that $\mc{R}_{z}$ lies within a proper special subvariety of $G$. By quotienting by the identity component of $\mc{R}$, we may assume that $\mc{R}$ is finite over $U$ and after finite base change that $\mc{R}$ is a section of $U$. Applying the Main Theorem of \cite{Mas89} and the extension to semi-abelian varieties given in Section $5$ of \emph{loc. cit.} to $\mc{R} \times \tau_1 \mc{R}, \dots, \tau_k \mc{R} \subset G^{k + 1}$, we get that $\mc{R}_z$ must lie within a proper special subvariety of $\mc{G}_z$.
\end{proof}

 \section{Applications of the main theorem}
 
In this section, we see how Theorem \ref{thm:main} can be used to get results on linear dependence of special points. We prove a weaker version of Corollary \ref{indpt1} which may be easier to apply. First we prove a lemma that is probably already known to experts but we could not find a reference for.

\begin{lem}\label{projection}
    Let $S$ be a Shimura variety and let $\pi_{ij} \colon S^n \to S^2$ be the projection onto the $i$th and $j$th coordinates. If $S' \subset S^n$ is a proper Shimura subvariety, then there exist $i, j \in \N$ with $1 \le i < j \le n$ such that $\pi_{ij}(S')$ is a proper Shimura subvariety of $S^2$.
\end{lem}
\begin{proof}
    We may assume that $S$ is a connected Shimura variety whose Shimura datum $(G_S, X)$ is such that $G_S$ is a product of simple groups and that $S'$ is also a connected Shimura variety whose Shimura datum is $(G', X')$. By viewing $S'$ as the orbit of a point under the action of a subgroup of $G_S^{n}$, it suffices to show that if $G' \subset G_1 \times \cdots G_n$ is a proper subgroup of a product of simple groups, then there exist $i, j$ such that $\pi_{ij}(G') \subset G_i \times G_j$ is a proper subgroup.

    Suppose otherwise, then since $\pi_{12}$ and $\pi_{13}$ are surjective, we can find elements of $G'$ of the form $(g, 1, g', \dots)$ and $(h, h', 1, \dots) \in G'$ for $g, h$ in $G_1, h' \in G_2, g' \in G_3$. Taking their commutator, we find that the projection of $G'$ to $G_1 \times G_2 \times G_3$ by taking the first three coordinates contains the element $(ghg^{-1}h^{-1}, 1, 1)$. Since $g, h \in G_1$ were arbitrary, we get that $[G_1, G_1] \times \set{1} \subset G_1 \times G_2 \times G_3$. Now since $G_1$ is a simple group, we get that $[G_1, G_1] = G_1$. By symmetry, we see that $\pi_{123}(G') = G_1 \times G_2 \times G_3$. By inductively repeating this process, we see that $G' = G_1 \times \cdots G_n$, contradicting the properness of $G'$.
\end{proof}

\begin{cor}\label{cor:LookAt2}
    Let $S, G, V$ be as in Corollary \ref{indpt1}. Fix $n \in \N$. There exists a proper Zariski closed subset $S' \subset S^2$ such that for any $n$ points $s_1, \dots, s_n \in S$ such that $(s_i, s_j) \not \in S'$ for any pair $i \neq j$, then any $V$-images $g_1, \dots, g_n$ for $s_1, \dots, s_n$ are linearly independent.
\end{cor}

In the case that $S = \mc{A}_2, \mc{A}_3$ is the moduli space of abelian surfaces or threefolds, we can give a slightly more explicit condition for the independence of $V$-images, by studying the subvarieties of $S^2$.

\begin{defn}\label{D independent}
    We say that a set of CM abelian varieties $s_1, \dots, s_n \in \mc{A}_g$ for $g \le 3$ are $D$-independent, for some integer $D$, if each $s_i$ is a simple abelian variety, there exists no isogeny of degree $\le D$ between some pair $s_i \to s_j$ for $i \neq j$, and for each $i$, we have $\Delta(s_i) > D$.
\end{defn}

We now give the proof of Corollary \ref{linearly independence of CM points}.

\begin{proof}[Proof of Corollary \ref{linearly independence of CM points}]
If $x_1, \cdots, x_n$ are linearly dependent in $G$ then by Corollary \ref{cor:LookAt2}, there must exist some $i \neq j$ such that $(s_i, s_j)$ lies within a finite set of proper special subvarieties of $S^2$. So we need to describe the special subvarieties in $\mc{A}_g \times \mc{A}_g$.

Let $W \subset \mc{A}_g \times \mc{A}_g$ be one of these finitely many proper special subvarieties. First, suppose that the projection to one of the factors is proper. This projection must be a proper special subvariety of $\mc{A}_g$. Since $g \le 3$, by \cite{MZ99}, the Shimura datum associated with the projection of $W$ must be PEL-type Shimura datum and so the endomorphism ring of some $s_i$, tensored up to $\Q$, must contain a $\Q$-algebra $B$. Since we assumed our $s_i$ were $D$-independent, they are simple and hence $B$ is a number field or division algebra. Since there are finitely many proper special subvarieties of $\mc{A}_g^2$, there are only finitely many such number fields or division algebras that can appear. Setting $D$ larger than the discriminants of the number fields and division algebras that appear prevent the projection of $s_i$ from being contained in a $W$ of this form.

Now suppose that the projection of $W$ to each factor of $\mc{A}_g$ is surjective. Choose a generic $w \in W$ and let $G \subset \Sp_{2g}(\C) \times \Sp_{2g}(\C)$ be the set of $g$ such that $gw \in W$. This is a Lie group of $\Sp_{2g}(\C) \times \Sp_{2g}(\C)$ and the projection onto each factor of $\Sp_{2g}(\C)$ is surjective. By Goursat's lemma, the group $G$ is given as the graph of an isomorphism between $\Sp_{2g}(\C)/N_1 \to \Sp_{2g}(\C)/N_2$, with $N_1, N_2$ normal subgroups. However, since $G$ is a Lie group and $\Sp_{2g}$ is a simple Lie group, we must have that $N_1, N_2$ are trivial and so $W$ is a Hecke correspondence. Therefore, by choosing $D$ large enough, we may also exclude all Hecke correspondences that appear by giving a lower bound on the isogenies that appear.

\end{proof}

 We give a generalization of \cite[Corollary 1.4]{PT22}. Although there are infinitely many abelian varieties with fixed discriminant of their generic $\mathbb{Z}$-Hodge structure, the number of irreducible special subvarieties that contains them is finite. 
 
\begin{proof}[Proof of Theorem \ref{Andre-Oort-Mordell-Lang}]
By applying Corollary \ref{cor:LookAt2} to $V^{r + 1}$ inside of $S^{r + 1} \times G^{r + 1}$, we obtain that there exists finitely many proper special subvarieties $S_n$ of $S^2$ such that for any $r + 1$ points $s_1, \dots, s_{r + 1} \in S$ with $(s_i, s_j) \not \in \bigcup_n S_n$ for all pairs $i, j$, any $V$-images of the $s_i$ are linearly independent. Since the Shimura datum is simple, the argument in the proof of Corollary \ref{linearly independence of CM points} gives that either the projection of $S_n$ on one factor is a proper special subvariety or $S_n$ is a correspondence on $S^2$. Let $S'$ be the union of all of the proper projections and let $D$ be the maximum degree of the projections of $S_n$ to a factor of $S$ when $S_n$ is a correspondence. Then, there exists an $N = N(r + 1, D)$ such that given $s_1, \dots, s_N \in S \bs S'$ special, we can find a subset $s_1', \dots, s_{r + 1}'$ such that no $(s_i', s_j')$ lies in a $S_n$, and thus their $V$-images are linearly independent. Thus, outside of $S' \subset S$, there are less than $N$ other special points of $S$ whose $V$-images lie in $\Gamma'$. Both $S'$ and $N$ depend on $S, G, V, r$, and are independent of $\Gamma$, proving the theorem.
\end{proof}

\subsection*{Acknowledgements} The authors were introduced to the problem at the 2023 Arizona Winter School under Jacob Tsimerman's project group. We would like to heartily thank Jacob for suggesting this question to us and for the helpful conversations, and Niven Achenjang, Camille Amoyal, Pam (Miao) Gu and Dong Gyu Lim for the helpful discussions during the winter school. We would also like to thank Sebastian Eterovi\'c for pointing out the application to multiplicative independence of singular moduli. Finally, we thank the organizers of the Arizona Winter School for creating a wonderful research experience.
\bibliographystyle{alpha}
\bibliography{ref}
\end{document}